# DISCUSSION OF "ANALYSIS OF VARIANCE—WHY IT IS MORE IMPORTANT THAN EVER" BY A. GELMAN

By Alan M. Zaslavsky

*Harvard University*

Andrew Gelman's contribution shifts the focus of "Analysis of Variance" (ANOVA) from the limited sense in which it has been commonly used in classical statistics, as a method of testing, to the broader framework of estimation and inference. The term more commonly used in this sense, "variance components modeling," also captures the same spirit. The essential idea is that of constructing distributions by using ideas of exchangeability; each variance component corresponds to a collection of exchangeable effects. This extremely powerful approach to linking the scientific structure of a dataset with a model has been and will continue to be widely applicable.

The transition from ANOVA to variance components modeling shifts attention from decomposition of the variance of the *sample* as in classical ANOVA to decomposition of the variance of the *population*. This shift in focus is appropriate to a world in which scientific questions become increasingly complex and are less frequently answerable through simple designed experiments.

As Gelman notes, in a Bayesian framework the estimation of variance components is relatively automatic; attention can be focused on defining sensible models rather than on constructing designs that can be analyzed easily.

Finally, Gelman's discussion of the manifold definitions of "fixed and random effects" is itself worth the price of admission.

**1. Testing and prior distributions.** Gelman reserves the issue of model specification, specifically testing of variance components, to a short section (Section 8.2) toward the end of the article (alluding to it only briefly in Section 3.4). Testing variance components is inherently different from testing a regression coefficient because the null hypothesis $\sigma^2 = 0$ is on the boundary of the parameter space. The difficulties this causes for hypothesis testing in the likelihood framework are well known.

In a Bayesian setting, we might distinguish two purposes of hypothesis testing: determining whether a scientifically interesting conclusion can be drawn with adequate certainty to be worth reporting, and selecting models (omitting unneeded model effects). For "scientific" testing of a regression coefficient we might select a locally uniform prior and then see whether we can at least be adequately confident of its sign, that is, is either $P(\beta > 0)$ or $P(\beta < 0)$ a posteriori close enough to 1? For a variance component, the boundary problem prevents defining a locally uniform prior or applying this "two-sided" approach. (Scale-invariant improper priors typically yield degenerate posteriors in variance components models.)

A "model-mixing" approach combines a point mass at the null with a proper distribution over the remainder of the distribution. I find this unsatisfactory as a default solution, especially in the context of independent priors on the magnitudes of variance components, because it requires informative prior beliefs about both the probability of the null hypotheses $\sigma_m^2 = 0$ (and the various combinations of nulls) and the scale of the variance component if nonzero.

We might avoid the scaling problem by defining prior distributions for *relative* variances $\sigma_m^2 / \sum_{m'} \sigma_{m'}^2$ rather than *absolute* variances $\sigma_m^2$. It is more natural to combine this prior with point masses for submodels because the prior probabilities for the submodels are relative to the distribution of a variable that is always scaled on $(0,1)$; the notion of a "small enough to be scientifically uninteresting" component is also more readily interpretable on a relative scale. Such a prior could also accommodate prior information about relative magnitudes of variance components as suggested by Gelman. Note that Gelman's suggestion of independent uniform priors on each component implies a uniform prior on these relative variance components, that is, a Dirichlet$(\delta, \delta, \ldots, \delta)$ prior with $\delta = 1$, conditional on the sum of the variance components (the marginal variance of the data in an additive model). A prior belief that the variance components should be nearly equal suggests a similar prior with $\delta \gg 1$, and if we believe that a few components should predominate, then we might assume $0 < \delta < 1$. At least in additive models, this prior specification allows us to separate specification of the prior for the marginal variance of the data from that for the ratios of components; such a separation is more difficult with independent priors for the different components.

I find posterior predictive tests [Rubin (1984) and Meng (1994)] a more satisfactory way to test variance components than model mixing: we fix a component at zero, and then by simulating data from its predictive distribution determine whether the observed value of a statistic related to that component in some monotone way is consistent with the predictions of the constrained model. Indeed, the sums-of-squares statistics of the classical ANOVA table are suitable for such a test. The boundary problem is not an issue with this approach since the reference distribution is determined by simulation rather than by asymptotics, and indeed no prior distribution is required for the variance component being tested. I would conjecture that for balanced data, the sums-of-squares statistics are optimal for posterior predictive testing of the null hypothesis on the corresponding variance components. An interesting research direction would be to prove this conjecture or find a superior statistic for the balanced case, and then to identify better statistics for posterior predictive testing with unbalanced data.

We might also be interested in testing as a means of selecting a model with fewer nuisance parameters, specifically by reducing the number of variance components. For this objective a conservative approach would incline toward retaining as many components as possible, but with a prior distribution that allows their estimates to stay close to zero if there is little or no evidence for nonzero values, possibly by using a small value of $\delta$ in a "relative" Dirichlet prior. Sensitivity of inferences of interest to the choice of prior might indicate that the data cannot unequivocally answer the questions of interest.

**2. Variance components as a focus of scientific research.** Much of the applied multilevel modeling literature treats variance components as nuisance parameters, putting the primary emphasis on estimation and testing of regression coefficients (representing scientifically interesting systematic relationships) or of functions of random effects (small area estimation in survey sampling and official statistics, profiling in health care, "league tables" for schools). An important exception is genetics, in which variance components are the basis for calculations of heritability. In my own applied research, I have found that variance components are also an inherently interesting object of inference. Two examples follow.

An analysis of predictors of administration of clinically appropriate chemotherapy for colorectal cancer estimated a residual variance component for hospital effects, after controlling for measured hospital and patient characteristics [Ayanian et al. (2003)]. To explain the importance of this variation to clinical readers, we noted that the difference between a moderately above-average and a moderately below-average hospital (1 SD above or below average) was about as large as the effect of the most important patient characteristic identified in the model. The large magnitude of this residual variation suggests that measurement of additional hospital characteristics might yield a scientific payoff. Furthermore, quality improvement activities might be directed to bringing lower-performing hospitals closer to the practices of their better-performing peers, consistent with arguments that substantial unexplained variation in rates of use of a medical procedure is in itself evidence of poor quality [Wennberg and Gittelsohn (1982)].

Samples of members of Medicare managed care health plans (private organizations that contract with the U.S. government to provide health care to elderly or disabled individuals) have been administered a survey annually for the last eight years to assess various aspects of the services they receive [Zaslavsky, Zaborski and Cleary (2004)]. The effects of measured characteristics of individual members or plans are fairly small, and the effects of an individual's characteristics are of little interest because the primary objective of the survey is to evaluate health care systems, not the predictors of an individual's reported experiences. These data were modeled with variance components for three levels of nested geographical units (region, state, Metropolitan Statistical Area or MSA) and for the organizational unit (the health plan). For ratings of "the plan" (primarily reflecting the quality of customer service interactions), the majority of variance (excluding the large bottom-level individual component) was explained by the organizational unit. However, the explainable variance for ratings of doctors was mainly attributable to geographical variation, with a smaller component attributable to the health plan. We interpreted this finding as reflecting the fact that the health plans have more control over customer services provided directly by the plan than over health care. The latter is largely provided by doctors and hospitals that are organizationally independent of the plans and might contract with multiple plans; other studies have shown substantial geographical variation in their practice patterns. This finding has implications for quality improvement, suggesting that interventions to improve quality of care might have to be directed to health care providers in an area rather than trying to identify and improve lower-performing health plans. Similar patterns were identified for other quality dimensions measured in the survey. A further analysis estimated variance components for the geographical units, the plan organization, time (year of survey administration) and interactions of these effects. The time effects were interesting in evaluating the extent to which relative changes in quality might be detected between consecutive years, while the plan by MSA interaction was useful for deciding whether to generate separate estimates by geographical area within large plans serving extensive areas. Estimation of these complex models was made possible by the unusual size of the survey datasets (over 700,000 respondents).

These examples illustrate that despite their relative unfamiliarity in many fields, variance components can be interpreted to nonstatisticians in a scientifically meaningful way.

## 3. Miscellaneous comments.

*Finite population variance components.* Gelman correctly notes that inference can be made for both finite-population and superpopulation variances, and that the distinction between these two targets of inference corresponds to the distinction between fixed and random effects. In my experience, estimation of a finite-population variance is relatively rarely of interest.

When we are really concerned about a specific set of units (such as alternative treatments in an experiment), we are likely to estimate rankings of and differences among those units as a basis for future action. On the other hand, when estimation of the variance component is intended as part of an inference about a more general law (as is common in econometric analyses of U.S. state data), we are likely to think of the finite population as part of a larger hypothetical population even if (as the 50 states) they in fact constitute the entire population.

*Method of moments.* Gelman draws out the connection between classical ANOVA and method-of-moments variance components estimators. Because these estimators are essentially linear combinations of variance statistics that can be directly calculated from the data, they have substantial heuristic value, since maximum-likelihood or Bayesian estimation in complex problems can be too much of a "black box" to yield adequate direct insight into the connection between the data and the parameter estimates. On the other hand, the simple decomposition of Gelman's equation (1) only applies when there are unbiased estimators with independent errors, which is not likely to be the case for complex models with unbalanced datasets.

DEPARTMENT OF HEALTH CARE POLICY
HARVARD UNIVERSITY MEDICAL SCHOOL
180 LONGWOOD AVENUE
BOSTON, MASSACHUSETTS 02115-5899
USA
E-MAIL: zaslavsk@hcp.med.harvard.edu